\documentclass[10pt,oneside,a4paper]{amsart}

\usepackage{amsmath,amsfonts,amssymb,amsthm}
\usepackage{array}
\usepackage[all,textures]{xy}
\usepackage{color}
\usepackage[orange]{xcolor}
\definecolor{blue(munsell)}{rgb}{0.0, 0.5, 0.69}
\usepackage{hyperref}
\hypersetup{hypertexnames = false, bookmarksdepth = 2, bookmarksopen = true, colorlinks, linkcolor = black, citecolor = blue(munsell), urlcolor = blue(munsell), pdfstartview={XYZ null null 1}}
\usepackage{tikz-cd}
\usepackage{etoolbox}
\usepackage{cleveref}
\usepackage{mathtools}
\usepackage{placeins}
\usepackage{caption}
\usepackage{enumitem}

\usepackage[T1]{fontenc}
\usepackage[utopia]{mathdesign}
\usepackage{eucal}

\usepackage{pdfsync}

\theoremstyle{plain}
\newtheorem{theorem}{Theorem}[section]
\newtheorem{proposition}[theorem]{Proposition}
\newtheorem{lemma}[theorem]{Lemma}

\theoremstyle{definition}
\newtheorem{definition}[theorem]{Definition}

\theoremstyle{remark}

\newtheorem{remark}[theorem]{Remark}

\newtheorem{notation}[theorem]{Notation}

\newcommand{\Mod}{\ensuremath{\mathsf{Mod}}}

\newcommand{\Hom}{\mathrm{Hom}}
\newcommand{\Obj}{\mathrm{Obj}}

\newcommand{\groth}{\ensuremath{\mathsf{Grt}}}
\newcommand{\site}{\ensuremath{\mathsf{Site}}}
\newcommand{\gsite}{\ensuremath{\mathsf{GrSite}}}
\newcommand{\topoi}{\ensuremath{\mathsf{Topoi}}}
\newcommand{\LC}{\ensuremath{\mathsf{LC}}}
\newcommand{\op}{\ensuremath{\mathsf{op}}}
\newcommand{\id}{\text{Id}}

\newcommand{\AAA}{\mathfrak{a}}
\newcommand{\BBB}{\mathfrak{b}}

\newcommand{\CCC}{\mathfrak{c}}
\newcommand{\DDD}{\mathfrak{d}}

\newcommand{\Sh}{\ensuremath{\mathsf{Sh}} }

\newcommand{\Fun}{\ensuremath{\mathsf{Fun}}}

\newcommand{\lra}{\longrightarrow}
\newcommand{\hlra}{\lhook\joinrel\longrightarrow}

\newcommand{\aaa}{\ensuremath{\mathcal{A}}}
\newcommand{\bbb}{\ensuremath{\mathcal{B}}}
\newcommand{\ccc}{\ensuremath{\mathcal{C}}}
\newcommand{\ddd}{\ensuremath{\mathcal{D}}}

\newcommand{\ttt}{\ensuremath{\mathcal{T}}}
\newcommand{\uuu}{\ensuremath{\mathcal{U}}}
\newcommand{\vvv}{\ensuremath{\mathcal{V}}}

\newcommand{\calR}{\ensuremath{\mathscr{R}}}

\newcommand{\calT}{\ensuremath{\mathscr{T}}}

\usetikzlibrary{arrows}
\tikzset{
	commutative diagrams/.cd,
	arrow style=tikz,
	diagrams={>=to}}

\AtBeginEnvironment{figure}{\stepcounter{equation}}
\AtBeginEnvironment{figure}{\stepcounter{figure}}

\CompileMatrices
\SelectTips{cm}{11}

\title{Grothendieck categories as a bilocalization of linear sites}
\author{Julia Ramos Gonz\'alez}
\address[Julia Ramos Gonz\'alez]{Universiteit Antwerpen, Departement Wiskunde-Informatica, Middelheimcampus,
	Middelheimlaan 1,
	2020 Antwerp, Belgium}
\email{julia.ramosgonzalez@uantwerpen.be}

\thanks{The author acknowledges the support of the Research Foundation Flanders (FWO) under Grant No. G.0112.13N}

\begin{document}
	\begin{abstract}
		We prove that the $2$-category $\groth$ of Grothendieck abelian categories with colimit preserving functors and natural transformations is a bicategory of fractions in the sense of Pronk of the $2$-category $\site$ of linear sites with continuous morphisms of sites and natural transformations.
		This result can potentially be used to make the tensor product of Grothendieck categories from earlier work by Lowen, Shoikhet and the author into a bi-monoidal structure on $\groth$.
	\end{abstract}	
	\maketitle
\section{Introduction}
Grothendieck categories are arguably the best-behaved and most studied large abel\-ian categories, second only to module categories which are their first examples. They play a fundamental role in algebraic geometry since the Grothendieck school, and are center stage in non-commutative algebraic geometry since the work of Artin, Stafford, Van den Bergh and others (see, for example, \cite{artintatevandenbergh}, \cite{artinzhang2}, \cite{staffordvandenbergh}). 
By the Gabriel-Popescu Theorem, Grothendieck categories can be viewed as ``linear topoi'', that is, as categories of sheaves on linear sites. The aim of this paper is to study the relation between Grothendieck categories and linear sites on a bicategorical level.

Throughout the paper $k$ will be a commutative ring. 

Let $\groth$ denote the 2-category of $k$-linear Grothendieck categories with colimit preserving $k$-linear functors and $k$-linear natural transformations. Let $\site$ denote the 2-category of $k$-linear sites with $k$-linear continuous morphisms of sites and $k$-linear natural transformations\footnote{We address size issues in \S \ref{section2catssitesgroth}.}. 

Given a linear site $(\AAA,\calT_{\AAA})$ we can naturally associate to it a Grothendieck category, namely its category of sheaves $\Sh(\AAA,\calT_{\AAA})$. In addition, given a continuous morphism $f: (\AAA,\calT_{\AAA}) \lra (\BBB, \calT_{\BBB})$ between sites, it naturally induces a colimit preserving functor $f^s: \Sh(\AAA,\calT_{\AAA}) \lra \Sh(\BBB, \calT_{\BBB})$ between the corresponding categories of sheaves (see \S \ref{seclinearsites}). In particular, among continuous morphisms there is a distinguished class of the so-called LC morphisms (see \Cref{defLC} below), which induce equivalences between the corresponding categories of sheaves. 

Observe that, from the Gabriel-Popescu theorem, it follows that every Grothendieck category can be realised as a category of sheaves on a linear site . Moreover, every colimit preserving functor between two categories of sheaves can be obtained as a ``roof'' of functors coming from continuous morphisms between sites, where the ``reversed arrows'' are equivalences induced by LC morphisms (see \Cref{roofthm} below).  These observations make it natural to view $\groth$ as a kind of ``localization'' of $\site$ at the class of LC morphisms. In this paper we make this idea precise by using the localization of bicategories with respect to a class of $1$-morphisms developed by Pronk in \cite{pronk} and further analysed by Tommasini in the series of papers \cite{tommasini1,tommasini2,tommasini3}. Our main result is the following:
\begin{theorem}[{\Cref{mainresultbis}}]\label{mainresult}
	There exists a pseudofunctor
	$$\Phi: \site \lra \groth$$
	which sends LC morphisms to equivalences in $\groth$, such that the pseudofunctor
	$$\tilde{\Phi}: \site[\LC^{-1}] \lra \groth$$
	induced by $\Phi$ via the universal property of the bicategory of fractions is an equivalence of bicategories.
\end{theorem}

The paper is structured as follows.

In sections \S \ref{seclinearsites} and \S \ref{sectionbicategoriesfractions} we provide the introductory material required for the rest of the paper.

Linear sites and Grothendieck categories are the linear counterpart of Grothendieck sites and Grothendieck topoi. In \S\ref{seclinearsites} we recall the basic notions of the corresponding linearized topos theory based on \cite{lowenlin} and \cite{lowen-ramos-shoikhet}. More concretely, we define continuous and cocontinuous morphism of linear sites and analyse the corresponding induced functors between the sheaf categories. In particular, we focus on the class of LC morphisms, as they play a fundamental role for the rest of the paper. 
 
Next, in \S\ref{sectionbicategoriesfractions}, we revisit some basic notions and results from \cite{pronk} and \cite{tommasini1} on bicategories and their localizations with respect to classes of 1-morphisms, which we will refer to as \emph{bilocalizations} from now on.

The tecnical core of the paper that allows the proof of \Cref{mainresult} is developed in \S\ref{section2catssitesgroth}, \S\ref{sectionlocalizingLC} and \S\ref{sectionlocalizationallsites}.
 
In \S\ref{section2catssitesgroth} we show that the natural map that assigns to each site its category of sheaves and to each continuous morphism between sites the induced colimit preserving functor between the categories of sheaves extends to a pseudofunctor
\begin{equation}
\Phi: \site \lra \groth.
\end{equation}
In addition, the class $\LC$ of LC morphisms admits a calculus of fractions in $\site$ and hence we have a bilocalization $\site[\LC^{-1}]$. This is done in \S\ref{sectionlocalizingLC}.

In particular, $\Phi$ sends LC morphisms to equivalences in $\groth$, and hence, by the universal property of the bilocalization we obtain a pseudofunctor from the bilocalization to $\groth$:
\begin{equation}
\tilde{\Phi}: \site[\LC^{-1}] \lra \groth.
\end{equation}

In \S\ref{sectionlocalizationallsites}, based on \cite{tommasini3}, we show  that the pseudofunctor $\Phi: \site \lra \groth$ fulfills the necessary and sufficient conditions for the induced pseudofunctor $\tilde{\Phi}: \site[\LC^{-1}] \lra \groth$ to be an equivalence of bicategories, which finishes the proof of \Cref{mainresult}. 

Our original interest in representing $\groth$ as a localization of $\site$ is of geometrical nature and comes from \cite{lowen-ramos-shoikhet}, where a tensor product of Grothendieck categories is defined in terms of a tensor product of linear sites. In particular, the fact that LC morphisms are closed under the tensor product of sites \cite[Prop 3.14]{lowen-ramos-shoikhet} is key to the proof of the independence of the tensor product of Grothendieck categories from the sheaf representations chosen. Combining \Cref{mainresult} with this fact makes it natural to think that the monoidal structure on $\site$ could be transferred to $\groth$ via the bilocalization, following the same principle as the \emph{monoidal localization} of ordinary categories from \cite{day}. This is briefly addressed in \S\ref{tensorproductvialocalization}. 

The reader can easily check that the methods and arguments used along the paper are also available in the classical set-theoretical setup of topos theory. Hence, if we denote by $\topoi$ the $2$-category of Grothendieck topoi with colimit preserving functors and natural transformations and by $\gsite$ the $2$-category of Grothendieck sites with continuous morphisms and natural transformations, we can state the corresponding analogue of \Cref{mainresult} above: 

\begin{theorem}
	There exists a pseudofunctor
	$$\Phi: \gsite \lra \topoi$$
	which sends LC morphisms to equivalences, such that the pseudofunctor
	$$\tilde{\Phi}: \gsite[\LC^{-1}] \lra \topoi$$
	induced by $\Phi$ via the universal property of the bicategory of fractions is an equivalence of bicategories.
\end{theorem}

\noindent \emph{Acknowledgement.} I am very grateful to Wendy Lowen for many interesting discussions, the careful reading of this manuscript and her valuable suggestions. I would also like to thank Ivo Dell'Ambrogio, Boris Shoikhet and Enrico Vitale for useful comments on bicategories of fractions and Matteo Tommasini, whose explanations on the behaviour of $2$-morphisms after localization have been essential in order to obtain the main result of this paper. 

\section{Linear sites and Grothendieck categories} \label{seclinearsites}
Linear sites and Grothendieck categories are the linear counterpart of the classical notions of Grothendieck sites and Grothendieck topoi from \cite{artingrothendieckverdier1}. We proceed now to give a brief account on this linearized version of topos theory. We refer the reader to \cite[\S 2]{lowenlin} for further details. 

Let $\AAA$ be a small $k$-linear category and $A \in  \AAA$ an object. A \emph{sieve on $A$} is a subobject $R$ of the representable module $\AAA(-,A)$ on $A$ in the category $\Mod(\AAA) \coloneqq \Fun_k(\AAA^{\op}, \Mod(k))$. 
Given $F = (f_i: A_i \rightarrow A)_{i \in I}$ a family of morphisms in $\AAA$, the \emph{sieve generated by $F$} is the smallest sieve $R$ on $A$ such that $f_i \in R(A_i)$ for all $i \in I$. We denote it by $\langle F \rangle = \langle f_i \rangle_{i \in I}$. 

A \emph{cover system} $\mathscr{R}$ on $\AAA$ consists of providing for each $A \in \AAA$ a family of sieves $\calR(A)$ on $A$. The sieves in a cover system $\calR$ are called \emph{covering sieves} or simply \emph{covers} (for $\calR$). We will say that a family $(f_i: A_i \lra A)_{i\in I}$ is a \emph{cover} if the sieve $\langle f_i \rangle_{i \in I}$ it generates is a cover.

A cover system $\calT$ on $\AAA$ is called a \emph{$k$-linear topology} if it fulfills the linearized version of the well-known identity, pullback and glueing axioms \cite[\S 2.2]{lowenlin}.
In particular, a \emph{$k$-linear site} is a pair $(\AAA,\calT)$ where $\AAA$ is a small $k$-linear category and $\calT$ is a $k$-linear topology on $\AAA$. 

As in the classical setting, one defines presheaves and sheaves as follows. 
\begin{definition}\label{defsheaf}
A \emph{presheaf} $F$ on $(\AAA, \calT)$ is simply an $\AAA$-module, i.e. $F$ is an object in $\Mod(\AAA)$.

A \emph{sheaf} $F$ on $(\AAA, \calT)$ is a presheaf such that the restriction functor
\begin{equation*}
F(A) \cong \Mod(\AAA)(\AAA(-,A), F) \lra \Mod(\AAA)(R, F)
\end{equation*}
is an isomorphism for all $A \in \AAA$ and all covering sieves $R \in \calT(A)$. 
We denote by $$\Sh(\AAA, \calT) \subseteq \Mod(\AAA)$$ the full subcategory of $k$-linear sheaves.
\end{definition}

In analogy with the classical setting, given a $k$-linear category $\AAA$, a $k$-linear topology $\calT$ on $\AAA$ is said to be \emph{subcanonical} if all the representable presheaves are sheaves for $\calT$. The finest $k$-linear topology on $\AAA$ for which all representable presheaves are sheaves is called the \emph{canonical topology} on $\AAA$.

We now proceed to give in more detail the corresponding linear notions of morphisms of sites and the morphisms of linear topoi induced by them. This is a linearized version of \cite[Expos\'e iii]{artingrothendieckverdier1}.

Given a $k$-linear functor $f: \AAA \lra \BBB$ between two $k$-linear categories $\AAA$ and $\BBB$, we have the following induced functors between their module categories:
\begin{itemize}
	\item $f^*: \Mod(\BBB) \lra \Mod(\AAA): F \longmapsto F \circ f$;
	\item Its left adjoint, denoted by $f_{!}: \Mod(\AAA) \lra \Mod(\BBB)$;
	\item Its right adjoint, denoted by $f_*: \Mod(\AAA) \lra \Mod(\BBB)$.
\end{itemize}	

\begin{definition}\label{defcontinuous}
	Consider $k$-linear sites $(\AAA,\calT_{\AAA})$ and $(\BBB,\calT_{\BBB})$. A $k$-linear functor $f: \AAA \lra \BBB$ is \emph{continuous}, if any of the following equivalent properties hold:
	\begin{enumerate}
		\item The functor $f^*:\Mod(\BBB) \lra \Mod(\AAA) : F \longmapsto F\circ f$ preserves sheaves;
		\item There exists a functor $f_s: \Sh(\BBB,\calT_{\BBB}) \lra \Sh(\AAA, \calT_{\AAA})$ such that the diagram
		\begin{equation*}\label{comsquare1}
		\begin{tikzcd}
		\Mod(\AAA) &\Mod(\BBB) \arrow[l,"f^*"']\\
		\Sh(\AAA, \calT_{\AAA}) \arrow[u,hook,"i_{\AAA}"] & \Sh(\BBB,\calT_{\BBB}) \arrow[l,"f_s"] \arrow[u,hook,"i_{\BBB}"']
		\end{tikzcd}
		\end{equation*}
		commutes;
		\item There exists a colimit preserving functor $f^s: \Sh(\AAA, \calT_{\AAA}) \lra \Sh(\BBB,\calT_{\BBB})$ such that the diagram
		\begin{equation*}\label{comsquare2}
		\begin{tikzcd}
		\AAA \arrow{r}{f} \arrow[d,hook,"Y_{\AAA}"'] &\BBB \arrow[d,hook,"Y_{\BBB}"]\\
		\Mod(\AAA) \arrow[r,"f_!"] \arrow[d,"\#_{\AAA}"'] &\Mod(\BBB) \arrow[d,"\#_{\BBB}"]\\
		\Sh(\AAA, \calT_{\AAA}) \arrow{r}{f^s} & \Sh(\AAA, \calT_{\AAA}) 
		\end{tikzcd}
		\end{equation*}
		commutes, where $Y_{\AAA}: \AAA \hlra \Mod(\AAA)$, $Y_{\BBB}: \BBB \hlra \Mod(\BBB)$ are the corresponding Yoneda embeddings and $\#_{\AAA}:\Mod(\AAA) \lra \Sh(\AAA,\calT_{\AAA})$, $\#_{\BBB}:\Mod(\BBB) \lra \Sh(\BBB,\calT_{\BBB})$ are the corresponding sheafification functors. 
	\end{enumerate}
	In addition, if any of the previous properties holds, we necessarily have that $f^s \dashv f_s$ and 
	\begin{equation}
	f^s \cong \#_{\BBB} \circ f_! \circ i_{\AAA}.
	\end{equation}
\end{definition}

\begin{definition}
Consider $k$-linear sites $(\AAA,\calT_{\AAA})$ and $(\BBB,\calT_{\BBB})$. A $k$-linear functor $f: \AAA \lra \BBB$ is \emph{cocontinuous}, if any of the following equivalent properties hold:
	\begin{enumerate}
		\item For each object $A \in \AAA$ and each covering sieve $R \in \calT_{\BBB} (f(A))$, there exists a covering sieve $S \in \calT_{\AAA}(A)$ with $f S \subseteq R$. 
		\item The functor $f_*: \Mod(\AAA, \calT_{\AAA}) \lra \Mod(\BBB, \calT_{\BBB})$ preserves sheaves.
	\end{enumerate}
\end{definition}	
In addition, if any of the previous properties holds we have:
\begin{enumerate}
	\item The functor $\widetilde{f}^* = \#_{\AAA} \circ f^* \circ i_{\BBB}$ is colimit preserving and exact and the diagram
	\begin{equation*}\label{comsquare3}
	\begin{tikzcd}
	\Mod(\AAA) \arrow[d,"\#_{\AAA}"'] &\Mod(\BBB) \arrow[d,"\#_{\BBB}"] \arrow[l,"f^*"']\\
	\Sh(\AAA, \calT_{\AAA}) &  \Sh(\BBB, \calT_{\BBB}) \arrow{l}{\widetilde{f}^*}
	\end{tikzcd}
	\end{equation*}
	is commutative up to canonical isomorphism. 
	\item There exists a functor $\widetilde{f}_*:\Sh(\AAA, \calT_{\AAA}) \lra \Sh(\BBB,\calT_{\BBB})$ such that the diagram
	\begin{equation*}\label{comsquare4}
	\begin{tikzcd}
	\Mod(\AAA, \calT_{\AAA})  \arrow{r}{f_*} &\Mod(\BBB, \calT_{\BBB}) \\
	\Sh(\AAA, \calT_{\AAA}) \arrow[r,"\widetilde{f}_*"'] \arrow[u,hook,"i_{\AAA}"] & \Sh(\BBB, \calT_{\BBB}) \arrow[u,hook,"i_{\BBB}"']
	\end{tikzcd}
	\end{equation*}
	commutes up to canonical isomorphism and $\widetilde{f}^* \dashv \widetilde{f}_*$ is an adjoint pair.
\end{enumerate}
\begin{remark}
	In \cite{lowenlin} the term \emph{cover continuous} is used for what we call here cocontinuous.
\end{remark}

Recall that a $k$-linear \emph{Grothendieck abelian category} $\ccc$ is a cocomplete abelian $k$-linear category with a generator and exact filtered colimits.

The fact that the categories of sheaves over linear sites (or linear Grothendieck topoi) are precisely the Grothendieck categories can be deduced from Gabriel-Popescu theorem \cite{gabrielpopescu} together with the main result in \cite{borceuxquinteiro}. Indeed, Gabriel-Popescu theorem characterizes Grothendieck categories as the localizations of presheaf categories of linear sites, that is subcategories of presheaf categories whose embedding functor has a left exact left adjoint. On the other hand, from \cite[Thm 1.5]{borceuxquinteiro} one deduces that categories of sheaves are precisely the localizations of presheaf categories of linear categories. Thus, the combination of the two results provides a linear counterpart of the classical Giraud Theorem that characterizes Grothendieck topoi in the classical setting. 

Observe, nevertheless, that the classical Gabriel-Popescu theorem does not provide us with all the possible realizations of Grothendieck categories as categories of linear sheaves. Such result is provided by the generalization of Gabriel-Popescu theorem in \cite{lowenGP}: given a Grothendieck category $\ccc$, it characterizes the linear functors $u:\AAA \lra \ccc$ such that the functor $$\CCC \lra \Mod(\AAA): C \longmapsto \ccc(u(-),C)$$ is a localization. 

We now introduce a distinguished class of continuous morphisms, called LC morphisms (see \cite[Def 3.4]{lowen-ramos-shoikhet}), where LC stands for ``Lemme de comparison'' (see \cite[\S 4]{lowenGP}). In particular, the functoriality of these morphisms with respect to the tensor product of linear sites constructed in \cite[\S 2.4]{lowen-ramos-shoikhet} is the key point in order to provide a well-defined tensor product of Grothendieck categories expressed in terms of realizations as sheaf categories \cite[\S 4.1]{lowen-ramos-shoikhet}. 

\begin{definition}\label{defLC}
	Consider a $k$-linear functor $f: \AAA \lra \CCC$. 
	\begin{enumerate}
		\item Suppose $\CCC$ is endowed with a cover system $\calT_{\CCC}$. We say that $f: \AAA \lra (\CCC, \calT_{\CCC})$ satisfies 
		\begin{itemize}
			\item[\textbf{(G)}] if for every $C \in \CCC$ there is a covering family $(f(A_i) \lra C)_i$ for $\calT_{\CCC}$.
		\end{itemize}
		\item Suppose $\AAA$ is endowed with a cover system $\calT_{\AAA}$. We say that $f: (\AAA, \calT_{\AAA}) \lra \CCC$ satisfies
		\begin{itemize}
			\item[\textbf{(F)}] if for every $c: f(A) \lra f(A')$ in $\CCC$ there exists a covering family $a_i: A_i \lra A$ for $\calT_{\AAA}$ and $f_i: A_i \lra A'$ with $cf(a_i) = f({f}_i)$;
			\item[\textbf{(FF)}] if for every $a: A \lra A'$ in $\AAA$ with $f(a) = 0$ there exists a covering family $a_i: A_i \lra A$ for $\calT_{\AAA}$ with $aa_i = 0$. 
		\end{itemize}
		\item Suppose $\AAA$ and $\CCC$ are endowed with cover systems $\calT_{\AAA}$ and $\calT_{\CCC}$ respectively. We say that $f: (\AAA, \calT_{\AAA}) \lra (\CCC, \calT_{\CCC})$ satisfies
		\begin{itemize}
			\item[\textbf{(LC)}] if $f$ satisfies (G) with respect to $\calT_{\CCC}$, (F) and (FF) with respect to $\calT_{\AAA}$, and we further have $\calT_{\AAA} = f^{-1} \calT_{\CCC}$.
		\end{itemize}
	\end{enumerate}
\end{definition}

The following gives a characterization of LC morphisms between linear sites in terms of continuity and cocontinuity, and it will be used in \S \ref{sectionlocalizingLC}. 
\begin{proposition}\label{LCcandcc}
	Consider a morphism of sites $f: (\AAA,\calT_{\AAA}) \lra (\BBB,\calT_{\BBB})$ satisfying \emph{\textbf{(G)}} with respect to $\calT_{\BBB}$ and \emph{\textbf{(F)}}, \emph{\textbf{(FF)}} with respect to $\calT_{\AAA}$. Then, the following are equivalent:
	\begin{enumerate}
		\item The morphism $f$ satisfies \emph{\textbf{(LC)}} (with respect to the topologies $\calT_{\AAA}$ and $\calT_{\BBB}$);
		\item The morphism $f$ is continuous and cocontinuous.
	\end{enumerate}
	\begin{proof}
	Assume (1) holds. The fact that $f$ is continuous is given by Lemme de comparison \cite[Cor 4.5]{lowenGP} and the fact that $f$ is cocontinuous follows from \cite[Lem 2.15]{lowenlin}.
	
	We prove now the converse. As by hypothesis $f$ is cocontinuous and satisfies \textbf{(F)} and \textbf{(FF)}, we have that $f^{-1}\calT_{\BBB} \subseteq \calT_{\AAA}$ by \cite[Prop 2.16]{lowenlin}. Now, as $f$ is continuous, by applying the linear counterpart of \cite[Expos\'e iii, Prop 1.6]{artingrothendieckverdier1}, we have that 
	$f(\calT_{\AAA}) \subseteq \calT_{\BBB}$. Consequently, by applying $f^{-1}$ we have that $\calT_{\AAA} \subseteq f^{-1}f \calT_{\AAA} \subseteq f^{-1}(\calT_{\BBB})$, which concludes the argument.
	\end{proof}
\end{proposition}

\begin{remark}\label{cocontandcontLC}
	Let $f: (\AAA,\calT_{\AAA}) \lra (\BBB,\calT_{\BBB})$ be an LC morphism. As it is continuous and cocontinuous, one can consider the induced functors between the corresponding sheaf categories both as a continuous and as a cocontinuous morphism. An easy check shows that those are related as follows:
	\begin{equation}
		\widetilde{f}_* \cong f^s,
	\end{equation}
	and hence 
	\begin{equation}
		\widetilde{f}^* \cong f_s.
	\end{equation}
\end{remark}

Our interest in LC morphisms is twofold. Firstly, they induce equivalences between the corresponding sheaf categories \cite[Cor 4.5]{lowenGP}. Secondly, making essential use of LC morphisms we are able recover any colimit preserving functor between Grothendieck categories as being induced by a roof of continuous morphisms of linear sites. More precisely:
\begin{theorem}[Roof theorem]\label{roofthm}
	Let $(\AAA, \calT_{\AAA})$ and $(\BBB, \calT_{\BBB})$ be linear sites and consider a colimit preserving functor $F: \Sh(\AAA,\calT_{\AAA}) \lra \Sh(\BBB,\calT_{\BBB})$. Then, there exist a subcanonical site $(\CCC,\calT_{\CCC})$ and a diagram
	\begin{equation}
	\begin{tikzcd}
	&\CCC\\
	\AAA \arrow[ur,"f"] &&\BBB, \arrow[ul,"w",swap]
	\end{tikzcd}
	\end{equation}
	where $f$ is a continuous morphism and $w$ is an LC morphism, such that
	\begin{equation}
	\begin{tikzcd}
	\Sh(\AAA, \calT_{\AAA}) \arrow[rr,"F"] \arrow[dr,"f^s",swap] &&\Sh(\BBB, \calT_{\BBB})\\
	&\Sh(\CCC,\calT_{\CCC}) \arrow[ur,"\tilde{w}^*",swap]
	\end{tikzcd}
	\end{equation}
	is a commutative diagram up to isomorphism. 
\end{theorem}

This theorem is a slight generalization of \cite[Tag 03A2]{stacksproject} in the linear setting, where the result is provided for geometric morphisms between Grothendieck topoi (i.e. adjunctions of functors between Grothendieck topoi where the left adjoint is left exact). Our version focuses on the left adjoints, or equivalently on the colimit preserving functors, without requiring them to be left exact (i.e. without requiring the adjunction to be a geometric morphism). Observe that we call LC morphism what in \cite{stacksproject} is called \emph{special cocontinuous functor} (see \Cref{LCcandcc}). 

The proof can be obtained along the lines of \cite[Tag 032A]{stacksproject}. We will just provide, for convenience of the reader, the construction of the site $(\CCC, \calT_{\CCC})$ and the morphisms $f$ and $w$, as these constructions will be frequently used throughout the paper.

Take $\CCC$ to be the full $k$-linear subcategory of $\Sh(\BBB,\calT_{\BBB})$ with the following set of objects
\begin{equation}
\Obj(\CCC)=\{\#_{\BBB}(\BBB(-,B))\}_{B\in \BBB} \cup \{F(\#_{\AAA}(\AAA(-,A)))\}_{A \in \AAA}.
\end{equation}
We endow it with the topology $\calT_{\CCC}$ induced from the canonical topology in $\Sh(\BBB,\calT_{\BBB})$.

Then we define $f: \AAA \lra \CCC$ as the composition $F \circ \#_{\AAA} \circ Y_{\AAA}$ and $w: \BBB \lra \CCC$ as the composition $\#_{\BBB} \circ Y_{\BBB}$.
 
 \section{Bicategories of fractions}\label{sectionbicategoriesfractions}
In this section we recall the main notions and results on localizations of bicategories from \cite{pronk} and \cite{tommasini1}. In general we will follow the notations and terminology from \cite{pronk} with the exception that, following a more standard terminology, we will call \emph{pseudofunctor} what in \cite{pronk} is called \emph{homomorphism of bicategories}.

We first fix some notations for the rest of the paper. 

Given a bicategory, we denote the vertical composition of $2$-morphisms by $\bullet$ and the horizontal composition of $2$-morphisms by $\circ$. In particular, given a diagram
\begin{equation*}
\begin{tikzcd}[column sep= 50pt ]
	A \arrow[r,bend left=40, "f"] \arrow[r,bend right=40,swap,"f"] \arrow[r, phantom,"\Downarrow \id_{f}", pos=0.6] &B  \arrow[r,bend left=40, "g"] \arrow[r,bend right=40, swap, "h"]  \arrow[r, phantom,"\Downarrow \alpha", pos=0.55] &C
\end{tikzcd}
\end{equation*}
in a bicategory $\ccc$, we denote by $\alpha \circ f$ to the horizontal composition $\alpha \circ \id_f$.

We will recall some important definitions for the rest of the paper.

\begin{definition}\label{defequiv}
	Given a 1-morphism $f: A \lra B$ in a bicategory $\ccc$, we say it is an \emph{equivalence} (or \emph{internal equivalence} in the terminology of \cite{tommasini1}) if there exists another 1-morphism $g: B \lra A$ and two invertible $2$-morphisms $\alpha: \id_A \Rightarrow g \circ f$ and $\beta: f \circ g \Rightarrow \id_B$ satisfying the triangle identities, i.e. the compositions
	\begin{equation*}
		\begin{tikzcd}[row sep=small]
			f \arrow[r, Rightarrow,"f \circ \alpha"] & f \circ g \circ f \arrow[r,Rightarrow, "\beta \circ f"] & f;\\
			g \arrow[r,Rightarrow,"\alpha \circ g"] & g \circ f \circ g \arrow[r,Rightarrow,"g \circ \beta"] & g
		\end{tikzcd}	
	\end{equation*}
	are the identity on $f$ and on $g$ respectively.
\end{definition}	

\begin{definition}\label{defpropertiespseudofunctor}
	A pseudofunctor $\Phi: \aaa \lra \bbb$ between two bicategories is
	\begin{itemize}
		\item \emph{essentially surjective on objects} if and only if for all $B \in \Obj(\bbb)$ there exists an $A \in\Obj(\aaa)$ such that there is an equivalence $\Phi(A) \cong B$ in $\bbb$; 
		\item \emph{essentially full} if for all $A, A' \in \Obj(\aaa)$ the functor 
		$$\Phi_{A,A'}: \aaa(A, A') \lra \bbb(\Phi(A), \Phi(A'))$$ 
		is essentially surjective ;
		\item \emph{fully faithful on 2-morphisms} if for all $A, A' \in \Obj(\aaa)$ the functor 
		$$\Phi_{A,A'}: \aaa(A, A') \lra \bbb(\Phi(A), \Phi(A'))$$ 
		is fully faithful;
		\item an \emph{equivalence of bicategories} if it is essentially surjective on objects, essentially full and fully faithful on 2-morphisms.
	\end{itemize}
\end{definition}	 

In \cite[\S 2]{pronk} a localization theory for bicategories along a class of 1-morphisms is developed generalizing the well-known localization of ($1$-)categories due to Gabriel-Zisman \cite{gabriel-zisman}. In particular, the bicategory of fractions in loc.cit. is defined and constructed by means of a right calculus of fractions. Observe that one could analogously develop the theory for a left calculus of fractions, as it is done in the $1$-categorical case. More precisely, a class of ($1$-)morphisms admits a left calculus of fractions if and only if the same class of ($1$-)morphisms in the opposite (bi)category admits a right calculus of fractions. Recall that the \emph{opposite bicategory} (or \emph{transpose bicategory} in the terminology of \cite{benabou}) is given by reversing the 1-morphisms and keeping the direction of the $2$-morphisms. In our case, we will use a left calculus of fractions, hence we introduce the analogous results from \cite{pronk} for a left calculus of fractions.
\begin{definition}{\cite[\S 2.1]{pronk}}
	Let $\ccc$ be a bicategory. We say a class $\mathsf{W}$ of 1-morphisms on $\ccc$ admits a \emph{left calculus of fractions} if it satisfies:
	\begin{itemize}
		\item[\textbf{LF1}] All equivalences belong to $\mathsf{W}$;
		\item[\textbf{LF2}] $\mathsf{W}$ is closed under composition of 1-morphisms;
		\item[\textbf{LF3}] Every solid diagram
		\begin{equation*}
			\begin{tikzcd}
				\AAA\\
				\CCC \arrow[u,"f"] \arrow[r,"w"'] & \BBB
			\end{tikzcd}
		\end{equation*}
		in $\ccc$ with $w \in \mathsf{W}$ can be completed to a square 
		\begin{equation*}
			\begin{tikzcd}
				\AAA \arrow[r,dotted,"v"] & \DDD\\
				\CCC \arrow[u,"f"] \arrow[r,"w"'] & \BBB \arrow[u,dotted,"g"'] \arrow[ul, start anchor={[xshift=-1.3ex, yshift=+1.3ex]}, end anchor={[xshift=+1.3ex, yshift=-1.3ex]}, Rightarrow,"\alpha"']
			\end{tikzcd}
		\end{equation*}
		where $\alpha$ is an invertible $2$-morphism and $v \in \mathsf{W}$;
		\item[\textbf{LF4}] \begin{itemize}
			\item[(1)] Given two morphisms $f,g : B \lra A$, a morphism $w: B' \lra B$ in $\mathsf{W}$ and a 2-morphism $\alpha: f \circ w \Rightarrow g \circ w$, there exists a morphism $v: A \lra A'$ in $\mathsf{W}$ and a 2-morphism $\beta: v \circ f \Rightarrow v \circ g$ such that $v \circ \alpha = \beta \circ w$;
			\item[(2)] if $\alpha$ is an isomorphism, we require $\beta$ to be an isomorphism too; and
			\item[(3)] given another pair $v': A \lra A'$ in $\mathsf{W}$ and $\beta': v' \circ f \Rightarrow v' \circ g$ satisfying condition (1), there exist 1-morphisms $u,u': A' \lra A''$ with $u \circ v$, $u' \circ v'$ in $\mathsf{W}$ and an invertible $2$-morphism $\epsilon: u \circ v \Rightarrow u' \circ v'$ such that the diagram
			\begin{equation*}
				\begin{tikzcd}
					u \circ v \circ f \arrow[r,"u \circ \beta"] \arrow[d,"\epsilon \circ f"'] & u \circ v \circ g \arrow[d,"\epsilon \circ g"]\\
					u' \circ v' \circ f \arrow[r,"u' \circ \beta'"'] &u' \circ v' \circ g
				\end{tikzcd}
			\end{equation*}
			is commutative;
		\end{itemize}
		\item[\textbf{LF5}] $\mathsf{W}$ is closed under invertible $2$-morphisms.
	\end{itemize}
\end{definition}

\begin{remark}
	The first axiom can be weakened as is done in \cite{tommasini1}.
\end{remark}

\begin{definition}{\cite[\S 2]{pronk}}
	Given a category $\ccc$ and a class of 1-morphisms $\mathsf{W}$ in $\ccc$ admitting a left calculus of fractions, a \emph{bilocalization of $\ccc$ along $\mathsf{W}$} is a pair $(\ccc[\mathsf{W}^{-1}],\Psi)$, where $\ccc[\mathsf{W}^{-1}]$ is a bicategory and $\Psi: \ccc \lra \ccc[\mathsf{W}^{-1}]$ is a pseudofunctor such that:
	\begin{enumerate}
		\item $\Psi$ sends elements in $\mathsf{W}$ to equivalences;
		\item Composition with $\Psi$ gives an equivalence of bicategories $$\Hom(\ccc[\mathsf{W}^{-1}], \ddd) \lra \Hom_{\mathsf{W}}(\ccc,\ddd)$$ for each bicategory $\ddd$, where $\Hom$ denotes the bicategory of pseudofunctors (see \cite[\S 8]{benabou}) and $\Hom_{\mathsf{W}}$ its full sub-bicategory of elements sending $\mathsf{W}$ to equivalences.
	\end{enumerate} 
	Observe that, in particular, $\ccc[\mathsf{W}^{-1}]$ is unique up to equivalence of bicategories \cite[\S 3.3]{pronk}.
\end{definition}

In \cite[\S2]{pronk} a detailed construction for $(\ccc[\mathsf{W}^{-1}],\Psi)$ is provided for a right calculus of fractions and in \cite{tommasini1} a simplified version of this construction is provided, less dependent of the axiom of choice. 
By inverting the direction of $1$-morphisms one gets the analogous construction of the bilocalization for a left calculus of fractions. 

 \section{The 2-category of Grothendieck categories and the 2-category of sites}\label{section2catssitesgroth} 

Fix a universe $\uuu$. For a $\uuu$-small $k$-linear site $(\AAA, \calT_{\AAA})$, the category $\Sh(\AAA, \calT_{\AAA})$ is defined with respect to the category $\uuu$-$\Mod(k)$ of $\uuu$-small $k$-modules. Let $\site$ denote the 2-category of $\uuu$-small $k$-linear sites with $k$-linear continuous morphisms of sites and $k$-linear natural transformations. By definition, a $\uuu$-Grothendieck category is a $k$-linear abelian category with a $\uuu$-small set of generators, $\uuu$-small colimits and exact $\uuu$-small filtered colimits. Let $\vvv$ be a larger universe such that all the categories $\Sh(\AAA, \calT_{\AAA})$ are $\vvv$-small and let $\groth$ denote the 2-category of $k$-linear $\vvv$-small $\uuu$-Grothendieck categories. Up to equivalence, $\groth$ is easily seen to be independent of the choice of $\vvv$. In the rest of the paper, we will omit the universes $\uuu$ and $\vvv$ from our notations and terminology.

\begin{remark}\label{enriched}
Observe that $\site$ and $\groth$ are actually enriched 2-categories, more precisely $k$-linear 2-categories in the sense of \cite[Def 2.4 \& 2.5]{ganter-kapranov}. 
\end{remark}

\begin{remark}
	Observe that equivalences (see \Cref{defequiv} above) in $\groth$ are just the colimit preserving $k$-linear functors which are equivalences of categories in the usual sense, while equivalences in $\site$ are just the $k$-linear continuous morphisms which are equivalences of categories in the usual sense.
\end{remark}

\begin{notation}
	We denote by $\LC$ the family of LC morphisms in $\site$ (see \Cref{defLC}). 
\end{notation}

The 2-category $\groth$ is related to $\site$ in a natural way. Indeed, we define a pseudofunctor
\begin{equation}\label{twofunctorsitegroth}
\Phi: \site \lra \groth
\end{equation}
as follows:
\begin{itemize}
	\item Given a site $(\AAA, \calT_{\AAA})$, we define
	$$\Phi(\AAA, \calT_{\AAA}) = \Sh(\AAA, \calT_{\AAA}),$$
	which is a Grothendieck category;
	\item Given a continuous map between two sites $f: (\AAA,\calT_{\AAA}) \lra (\BBB,\calT_{\BBB})$, we define 
	$$\Phi(f) = \begin{tikzcd}
	\Sh(\AAA, \calT_{\AAA}) \arrow[r,"f^s"] &\Sh(\BBB,\calT_{\BBB}),
	\end{tikzcd}$$
	which is colimit preserving;
	\item  Given two continuous morphisms $f,g : (\AAA,\calT_{\AAA}) \lra (\BBB, \calT_{\BBB})$ and a natural transformation $\alpha: f \Rightarrow g$, we define a natural transformation 
	\begin{equation}\label{thenaturaltransform}
		\Phi(\alpha)= \alpha^s: f^s \Rightarrow g^s
	\end{equation}
	as follows. For any $F \in \Sh(\BBB,\calT_{\BBB})$ and any $A \in \AAA$, we have the following morphism:
	$$(\alpha_s)_F (A) \coloneqq F(\alpha_A): g_s(F)(A) = F(g(A)) \lra F(f(A)) = f_s(F)(A)$$  
	which is $k$-linear and natural in $A$ and $F$ and hence it defines a natural transformation $\alpha_s: g_s \Rightarrow f_s$. We define $\alpha^s$ as the natural transformation corresponding to $\alpha_s$ via the natural adjunctions. More precisely, for all $F \in \Sh(\BBB,\calT_{\BBB})$ and all $G \in \Sh(\AAA, \calT_{\AAA})$ we have a composition:
	\begin{equation*}
	\begin{tikzcd}
	\Sh(\AAA,\calT_{\AAA})(G,g_s(F)) \arrow[r,"(\alpha_s)_F \circ \: -"]  &\Sh(\AAA,\calT_{\AAA})(G, f_s(F)) \arrow[d,"\cong"]\\
	\Sh(\BBB,\calT_{\BBB})(g^s(G),F) \arrow[u,"\cong"] &\Sh(\BBB,\calT_{\BBB})(f^s(G), F)
	\end{tikzcd}
	\end{equation*}
	where the vertical functors are the adjunctions. Observe this composition is natural in $F$ and $G$. Consequently, there is an induced 2-morphism $f^s \Rightarrow g^s$, and this is the 2-morphism we denote by $\alpha^s$.
\end{itemize} 
One can easily check these data indeed define a pseudofunctor:
\begin{itemize}
	\item Given any two sites $(\AAA,\calT_{\AAA}), (\BBB,\calT_{\BBB})$ the map
	$$\Phi_{\AAA,\BBB}: \site(\AAA, \BBB) \lra \groth(\Sh(\AAA), \Sh(\BBB))$$
	induced by $\Phi$ is a functor. Indeed, consider a continuous morphism $f: \AAA \lra \BBB$, trivially the 2-morphism $\id_{f}: f \Rightarrow f$ is mapped to $\id_{f^s}: f^s \Rightarrow f^s$. Now consider $\alpha: f \Rightarrow g$ and $\beta: g \Rightarrow h$ and their vertical composition $\beta \bullet \alpha:f \Rightarrow h$. One has that  $$\left[ \alpha_s \bullet \beta_s \right]_G(A) = G(\alpha_A) \circ G(\beta_A) = G((\alpha \bullet \beta)_A )= \left[ (\beta \bullet \alpha)_s\right]_G(A)$$ 
	for all $G \in \Sh(\BBB)$ and all $A \in \AAA$. Hence, by adjunction, $\Phi_{\AAA,\BBB}$ preserves compositions.
	\item Let $\AAA= (\AAA,\calT_{\AAA})$ be a site and consider its identity morphism $\id_{\AAA}$. One has that $(\id_{\AAA})_s = \id_{\Sh(\AAA)}$ is the identity functor of $\Sh(\AAA)$. Hence, by adjunction, $$\id_{\Sh(\AAA)} \cong (\id_{\AAA})^s,$$ which gives us the unitor of $\Phi$. 
	
	Consider now two continous morphisms $f:\AAA \lra \BBB$ and $g: \BBB \lra \CCC$ in $\site$. By definition, we have that $$(g\circ f)^s \cong g^s \circ f^s,$$ which provides the associator of $\Phi$.
	\item It can be readily seen, using the fact that adjoints are unique up to unique isomorphism, that the unitor and associator of $\Phi$ fulfill the corresponding coherence axioms. 
\end{itemize}

Observe that $\Phi$ sends LC morphisms to equivalences. This is a direct consequence of the Lemme de comparaison. Hence, if $\LC$ admits a left calculus of fractions in $\site$, we will get, by the universal property of bilocalizations, a pseudofunctor
\begin{equation}
\tilde{\Phi}: \site[\LC^{-1}] \lra \groth.	
\end{equation}

\begin{remark}
	Recall from \Cref{enriched} that $\site$ and $\groth$ are $k$-linear $2$-categories, and observe that $\Phi$ is also a $k$-linear pseudofunctor. While in this paper we only need the bilocalization to exist as an ordinary bicategory, it is possible to show that in this case the bilocalization automatically satisfies the universal property of an ``enriched bilocalization'' (and in particular, the induced functor $\tilde{\Phi}$ is automatically $k$-linear), where we use the term in analogy with the enriched localizations from \cite{wolff}.
\end{remark}

\section{Bilocalization of  the 2-category of sites with respect to LC morphisms}\label{sectionlocalizingLC}
In this section we prove that $\LC$ admits a left calculus of fractions in $\site$.

First, we fix the following notations. Given a site $(\AAA, \calT_{\AAA})$, it will usually be denoted simply by $\AAA$ for the sake of brevity. We will denote by $i_{\AAA}: \Sh(\AAA) \hlra \Mod(\AAA)$ the natural inclusion and by $\#_{\AAA}:\Mod(\AAA) \lra \Sh(\AAA)$ the corresponding sheafification functor. The indexes will be omitted if the site we are working with is clear from the context. Furthermore, given an object $A \in \AAA$ we will denote $h_A = \AAA(-,A)$ the corresponding representable presheaf and by $h^{\#}_A = \#(\AAA(-,A)) $ its sheafification. We adopt these notations in order to simplify the formulas that will appear further in the paper. 

\begin{lemma}\label{CondLF1}
	Condition \emph{\textbf{LF1}} holds for $\LC$ in $\site$.
	\begin{proof}
		Take an equivalence $f \in \site(\AAA, \BBB)$, and denote by $g: \BBB \lra \AAA$ its quasi-inverse. Then, it is easy to see that the induced functor 
		$$f_s: \Sh(\BBB)  \lra \Sh(\AAA)$$
		is an equivalence of Grothendieck categories, with quasi-inverse given by $g_s$. We prove that $f$ belongs to $\LC$. Property \textbf{(G)} follows from the fact that $f$ is essentially surjective, and properties \textbf{(F)} and \textbf{(FF)} follow immediately from the fact that $f$ is fully-faithful. By \Cref{LCcandcc}, it only remains to prove that $f$ is cocontinuous. One can easily see that $$f^*: \Mod(\BBB) \lra \Mod(\AAA)$$ is also an equivalence with quasi-inverse given by $g^*: \Mod(\AAA) \lra \Mod(\BBB)$. Hence we have that $f_* \cong f_! \cong g^*$ by unicity of adjoint functors. Observe then that 
		$$f_* \circ i_{\AAA} \cong g^* \circ i_{\AAA} = i_{\BBB} \circ g_s,$$
		which implies that $f_*|_{\Sh(\AAA)}$ takes values in $\Sh(\BBB) \subseteq \Mod(\BBB)$. Consequently $f$ is cocontinuous. 		
	\end{proof}
\end{lemma}

\begin{lemma}\label{CondLF2}
	Condition \emph{\textbf{LF2}} holds for $\LC$ in $\site$.
	\begin{proof}
		The composition of continuous morphisms is again continuous, and the analogous statement is true for cocontinuous morphisms. So we only have to see that the composition of two LC morphisms again fulfills properties \textbf{(G)}, \textbf{(F)} and \textbf{(FF)} and we conclude by \Cref{LCcandcc}.
		
		Consider $v: \AAA \lra \BBB$ and $w: \BBB \lra \CCC$ two LC morphisms between sites.
		
		Property \textbf{(G)} for $w \circ v$ follows immediately from the fact that both $v$ and $w$ have property \textbf{(G)} and that the composition of coverings is again a covering (this is a direct consequence of the glueing axiom \cite[\S 2.2 ]{lowenlin}). 
		
		We now prove property \textbf{(F)}. Consider a morphism 
		$$c: (w \circ v) (A) \lra (w \circ v) (A')$$
		in $\CCC$. As $w$ has property \textbf{(F)}, we know there exist a covering $\{r_i: B_i \lra v(A)\}_i$ in $\BBB$ and morphisms $b_i:B_i \lra v(A')$ such that 
		\begin{equation}\label{propF}
		c \circ w(r_i) = w(b_i).
		\end{equation}
		Now, as $v$ has property \textbf{(G)}, we can find covering families $\{s_{ij}:v(A_{ij}) \lra B_i\}_j$ in $\BBB$ for all $i$. Now consider the morphisms $r_i \circ s_{ij}: v(A_{ij}) \lra v(A)$ and $b_i \circ s_{ij}:v(A_{ij}) \lra v(A')$ for each $i,j$. As $v$ has property \textbf{(F)}, we know there exist coverings $\{t_{ijk} : A_{ijk} \lra A_{ij}\}_{k}$ and $\{t'_{ijk}: A'_{ijk} \lra A_{ij}\}_{k}$ and morphisms $a_{ijk}: A_{ijk} \lra A$ and $a'_{ijk}: A'_{ijk} \lra A'$ in $\AAA$, such that 
		\begin{equation*}
		\begin{aligned}
		r_i \circ s_{ij} \circ v(t_{ijk}) &= v(a_{ijk})\\
		b_i \circ s_{ij} \circ v(t'_{ijk}) &= v(a'_{ijk})
		\end{aligned}
		\end{equation*} 
		for all $i,j,k$. Consider now the intersection of the two covering sieves generated by $\{t_{ijk} : A_{ijk} \lra A_{ij}\}_{k}$ and $\{t'_{ijk}: A'_{ijk} \lra A_{ij}\}_{k}$ in $\AAA$ for each $i,j$ and take a covering family $\{u_{ijk}:A_{ijk} \lra A_{ij}\}_k$ generating this sieve. Then we have that
		\begin{equation*}
		\begin{aligned}
		r_i \circ s_{ij} \circ v(u_{ijk}) &= v(\bar{a}_{ijk})\\
		b_i \circ s_{ij} \circ v(u_{ijk}) &= v(\bar{\bar{a}}_{ijk})
		\end{aligned}
		\end{equation*} 
		for morphisms $\bar{a}_{ijk}: A_{ijk} \lra A$ and $\bar{\bar{a}}_{ijk}: A_{ijk} \lra A'$. Observe that the family formed by the compositions $\{r_i \circ s_{ij} \circ v(u_{ijk})\}_{i,j,k}$ is a covering because it is a composition of coverings (note that $\{v(u_{ijk})\}_k$ is a covering of $v(A_{ij})$ in $\BBB$ because $\{u_{ijk}\}_k$ is a covering of $A_{ij}$ in $\AAA$ and $v$ is an LC morphism). Consequently $\{\bar{a}_{ijk}\}_{i,j,k}$ is a covering in $\AAA$ of $A$, because $v$ is an LC morphism and $\{v(\bar{a}_{ijk})\}_{i,j,k}$ is a covering on $\BBB$. Hence precomposing with $w(s_{ij} \circ v(u_{ijk}))$ in both terms of (\ref{propF}) we have that:
		$$c \circ w(r_i) \circ w(s_{ij} \circ v(u_{ijk})) = w(b_i) \circ w(s_{ij} \circ v(u_{ijk}))$$ 
		for all $i,j,k$. Observe that the first term is equal to $c \circ (w\circ v)(\bar{a}_{ijk})$ and the second term is equal to $(w \circ v)(\bar{\bar{a}}_{ijk})$. This proves \textbf{(F)} for $w \circ v$.
		
		To conclude, we prove property \textbf{(FF)}. Consider a morphism $a:A \lra A'$ in $\AAA$ such that $(w \circ v)(a) = 0$. As $w$ has property \textbf{(FF)}, there is a covering $\{r_i:B_i \lra v(A)\}_i$ in $\BBB$ such that
		\begin{equation}\label{propFF}
		v(a) \circ r_i = 0 
		\end{equation}
		for all $i$. As $v$  has property \textbf{(G)}, for each $i$ there exists a covering $\{s_{ij}: v(A_{ij}) \lra B_i\}_j$ in $\BBB$. Consider the covering $\{r_i \circ s_{ij}: v(A_{ij}) \lra v(A)\}_{i,j}$ of $v(A)$ in $\BBB$ given by the composition. As $v$ has property \textbf{(F)}, for each $i,j$ there exist a covering $\{t_{ijk}:A_{ijk} \lra A_{ij}\}_k$ in $\AAA$ and a family of morphisms $\bar{a}_{ijk}: A_{ijk} \lra A$ such that:
		\begin{equation} \label{confirmcovering}
		r_i \circ s_{ij} \circ v(t_{ijk}) = v(\bar{a}_{ijk}).
		\end{equation}
		Then, precomposing in both terms of (\ref{propFF}) with $s_{ij} \circ v(t_{ijk})$, one has that
		\begin{equation*}
		v(a) \circ r_i \circ s_{ij} \circ v(t_{ijk}) = v(a \circ \bar{a}_{ijk})= 0.
		\end{equation*}
		Eventually, as $v$ has property \textbf{(FF)}, we know that for every $i,j,k$ there exists a covering $\{u_{ijkl}:A_{ijkl} \lra A_{ijk}\}_{l}$ in $\AAA$ such that: 
		\begin{equation*}
		a \circ \bar{a}_{ijk} \circ u_{ijkl} = 0
		\end{equation*}
		But the family $\{\bar{a}_{ijk} \circ u_{ijkl}:A_{ijkl} \lra A\}_{ijkl}$ is a covering because it is a composition of coverings (observe in (\ref{confirmcovering}) that, as $v$ is LC, $\{v(t_{ijk})\}_{k}$ is a covering in $\BBB$ and hence so is $\{v(\bar{a}_{ijk})\}_{i,j,k}$ and thus $\{\bar{a}_{ijk}\}_{i,j,k}$ is a covering in $\AAA$). Hence we conclude the argument.
	\end{proof}
\end{lemma}

\begin{lemma}\label{CondLF3}
	Condition \emph{\textbf{LF3}} holds for $\LC$ in $\site$.
	\begin{proof}
		Assume we have a solid diagram 
		\begin{equation*}
		\begin{tikzcd}
		\BBB\\
		\CCC \arrow[u,"w"] \arrow[r,"f"'] & \AAA
		\end{tikzcd}
		\end{equation*}
		in $\site$ with $w \in \LC$. We have to prove that it can be completed to a square 
		\begin{equation*}
		\begin{tikzcd}
		\BBB \arrow[r,dotted,"g"] \arrow[dr, start anchor={[xshift=+1.3ex, yshift=-1.3ex]}, end anchor={[xshift=-1.3ex, yshift=+1.3ex]}, Rightarrow,"\alpha"]& \DDD\\
		\CCC \arrow[u,"w"] \arrow[r,"f"'] & \AAA \arrow[u,dotted,"v"'] 
		\end{tikzcd}
		\end{equation*}
		where $\alpha$ is an invertible $2$-morphism and $v \in \LC$.
		
		Consider the following morphism of Grothendieck categories
		\begin{equation*}
		\begin{tikzcd}
		\Sh(\AAA) \arrow[r, "f_s"] &\Sh(\CCC) \arrow[r,"\tilde{w}_*"]  &\Sh(\BBB)
		\end{tikzcd}
		\end{equation*}
		induced by $f$ and $w$, whose left adjoint is given by $f^s \circ \tilde{w}^*: \Sh(\BBB) \lra \Sh(\AAA)$ (see \S \ref{seclinearsites}).
		
		We apply the roof theorem (\Cref{roofthm}) to this latter morphism. Take $\DDD$ the site defined as the full subcategory of $\Sh(\AAA)$ with objects $\{h^{\#}_A\}_{A \in \AAA} \cup \{(f^s \circ \tilde{w}^*)(h^{\#}_{B})\}_{B \in \BBB}$ and we endow it with the topology induced by the canonical topology in $\Sh(\AAA)$. Then consider the following morphisms:
		\begin{equation*}
		\begin{tikzcd}
		\BBB \arrow[r,"g"] & \DDD\\
		& \AAA \arrow[u,"v"']
		\end{tikzcd}
		\end{equation*}
		defined $v = \#_{\AAA} \circ Y_{\AAA}$ and $g = (f^s \circ \tilde{w}^*) \circ \#_{\BBB} \circ Y_{\BBB}$. 
		
		Now, as $w$ is an LC morphism by hypothesis, and thus continuous and cocontinuous, we have the following chain of invertible $2$-morphisms (see \S \ref{seclinearsites} for the properties of continuous and cocontinuous functors):
		\begin{equation*}
		\begin{aligned}
		g \circ w &= f^s \circ \tilde{w}^* \circ \#_{\BBB} \circ Y_{\BBB} \circ w \\
		&=  f^s \circ \#_{\CCC} \circ w^*\circ i_{\BBB} \circ \#_{\BBB} \circ Y_{\BBB} \circ w\\
		&= f^s \circ \#_{\CCC} \circ w^*\circ i_{\BBB} \circ w^s \circ \#_{\CCC} \circ Y_{\CCC}\\
		&= f^s \circ \#_{\CCC} \circ i_{\CCC} \circ w_s \circ w^s \circ \#_{\CCC} \circ Y_{\CCC}\\
		&\cong f^s \circ \#_{\CCC} \circ Y_{\CCC}\\  
		&= \#_{\AAA} \circ Y_{\AAA} \circ f\\
		&= v \circ f
		\end{aligned}
		\end{equation*}
		Hence, if we take $\alpha$ to be this invertible $2$-morphism, we conclude the argument.
	\end{proof}
\end{lemma}

\begin{lemma}\label{CondLF5}
	Condition \emph{\textbf{LF5}} holds for $\LC$ in $\site$.
	\begin{proof}
		Consider two morphisms $v,w: \AAA \lra \BBB$ in $\site$ such that $w \in \LC$. Assume we also have an invertible $2$-morphism $\alpha: v \Rightarrow w$. We want to prove that $v$ also belongs to $\LC$.
		
		Consider an object $B \in \BBB$. By hypothesis, there exists a covering $\{r_i: w(A_i) \lra B \}_i$ and we can consider the associated covering 
		\begin{equation*}
			\{v(A_i) \xrightarrow{\alpha_{A_i}} w(A_i) \xrightarrow{r_i} B \}
		\end{equation*} 
		by composing with the isomorphism $\alpha_{A_i}$ given by the invertible $2$-morphism (recall that any isomorphism generates a covering sieve in a topology, the corresponding representable one). This proves property \textbf{(G)} for $v$.
		
		Consider now a morphism $b: v(A) \lra v(A')$ in $\BBB$. Then we can take the morphism $\alpha_{A'} \circ b \circ \alpha_{A}^{-1}: w(A) \lra w(A')$, and as property \textbf{(F)} holds for $w$, we have that there exists a covering $\{s_i:A_i \lra A\}_i$ and morphisms $a_i: A_i \lra A'$ in $\AAA$, such that 
		\begin{equation*}
			\alpha_{A'} \circ b \circ \alpha_{A}^{-1} \circ w(s_i) = w(a_i)
		\end{equation*}
		for all $i$. Consequently:
		\begin{equation*}
			v(a_i) = \alpha_{A'}^{-1} \circ w(a_i) \circ \alpha_{A_i} = \alpha_{A'}^{-1} \circ \alpha_{A'} \circ b \circ \alpha_{A}^{-1} \circ w(s_i) \circ \alpha_{A_i} = b \circ v(s_i),
		\end{equation*} 
		which proves that property \textbf{(F)} holds for $v$.
		
		Take now $a:A \lra A'$ in $\AAA$ such that $v(a)= 0$. This implies that $w(a) = 0$ and hence there exists a covering $\{t_i:A_i \lra A \}_i$ such that $a \circ t_i = 0$ for all $i$, which proves that property \textbf{(FF)} holds for $v$. 
		
		It remains to prove that $\calT_{\AAA} = v^{-1}\calT_{\BBB}$. We have that $\calT_{\AAA} = w^{-1} \calT_{\BBB}$ by hypothesis, hence given a covering sieve $\langle r_i: A_i \lra A \rangle$ of $A$ in $\AAA$, we have that $\langle r_i \rangle \in \calT_{\AAA}(A)$ if and only if $\langle w(r_i): w(A_i) \lra w(A) \rangle \in \calT_{\BBB}(w(A))$. 
		
		On the other hand, it is easy to see that $\langle w(r_i) \rangle \in \calT_{\BBB}(w(A))$ if and only if $$\langle \alpha_A^{-1} \circ w(r_i) \circ \alpha_{A_i} = v(r_i): v(A_i) \lra v(A) \rangle \in \calT_{\BBB}(v(A)),$$ because $\alpha_{A_i}, \alpha_A$ are isomorphisms. This concludes the argument.
	\end{proof}
\end{lemma}

\begin{lemma}\label{CondLF4}
	Condition \emph{\textbf{LF4}} holds for $\LC$ in $\site$.
	\begin{proof}
		First we prove (1) holds. Given two continuous morphisms $f,g : \BBB \lra \AAA$ in $\site$, an LC morphism $w: \BBB' \lra \BBB$ and a 2-morphism $\alpha: f \circ w \Rightarrow g \circ w$, we have to prove that there exists an LC morphism $v: \AAA \lra \AAA'$ and a 2-morphism $\beta: v \circ f \Rightarrow v \circ g$ such that $v \circ \alpha = \beta \circ w$.
		
		Consider the morphisms
		\begin{equation*}
		\begin{tikzcd}
		\Sh(\AAA) \arrow[r,"(f \circ w)_s"] &\Sh(\BBB')  &\text{and} &\Sh(\AAA) \arrow[r,"(g \circ w)_s"] &\Sh(\BBB') 
		\end{tikzcd}
		\end{equation*}
		between the corresponding sheaf categories, with left adjoints $(f \circ w)^s: \Sh(\BBB') \lra \Sh(\AAA)$ and $(g \circ w)^s:\Sh(\BBB') \lra \Sh(\AAA)$ respectively. We now perform a similar construction to that on the roof theorem. Take the full subcategory $\AAA'$ of $\Sh(\AAA)$ with $$\Obj(\AAA') = \{h^{\#}_{A}\}_{A \in \AAA} \cup \{(f \circ w)^s (h^{\#}_{B'})\}_{B' \in \BBB'} \cup \{(g \circ w)^s (h^{\#}_{B'})\}_{B' \in \BBB'},$$ and endow it with the topology given by the restriction of the canonical topology in $\Sh(\AAA)$. In particular, $\AAA'$ is subcanonical.
		We construct the following roofs
		\begin{equation*}
		\begin{tikzcd}
		&\AAA' &&&\AAA'\\
		\BBB' \arrow[ur,"r_f"] &&\AAA \arrow[ul,"v"'] &\BBB' \arrow[ur,"r_g"] &&\AAA \arrow[ul,"v"']
		\end{tikzcd}
		\end{equation*}
		where $v = \#_{\AAA} \circ Y_{\AAA}$, $r_f =(f \circ w)^s \circ \#_{\BBB'} \circ Y_{\BBB'}$ and $r_g= (g \circ w)^s \circ \#_{\BBB'} \circ Y_{\BBB'}$. One can easily see, following the same arguments as in the roof theorem, that $v$ is an LC morphism and that $(f \circ w)^s \cong \tilde{v}^* \circ (r_f)^s$ and $(g \circ w)^s \cong \tilde{v}^* \circ (r_g)^s$. We have chosen this special $\AAA'$ in order to have the same site on the top of both roofs, but the reasoning to prove that these roofs behave as the usual roof construction is not affected by this enlargement of the top category, as it remains to be small. Now that $v$ is constructed, we proceed to build the 2-morphism $\beta: v \circ f \Rightarrow v \circ g$.
		
		First observe that given $\alpha: f \circ w \Rightarrow g \circ w$ we have  $\alpha^s: (f \circ w)^s \Rightarrow (g \circ w)^s$ the induced 2-morphism described in \eqref{thenaturaltransform}. In particular, one has that:
		\begin{equation}\label{alphasrepres}
			\begin{tikzcd}
				\Hom_{\Sh(\BBB')}(h^{\#}_{(g \circ w)(A)},F) \arrow[r,"- \: \circ h^{\#}_{\alpha_A}"] \arrow[d,"\cong"] &\Hom_{\Sh(\BBB')}(h^{\#}_{(f\circ w)(A)}, F) \arrow[d,"\cong"]\\
				F((g \circ w)(A)) \arrow[r,"F(\alpha_A)"] \arrow[d,"\cong"] &F((f \circ w)(A)) \arrow[d,"\cong"]\\
				\Hom_{\Sh(\AAA)}(h^{\#}_A,(g \circ w)_s(F)) \arrow[r,"(\alpha_s)_F \circ \: -"] \arrow[d,"\cong"] &\Hom_{\Sh(\AAA)}(h^{\#}_A, (f\circ w)_s(F)) \arrow[d,"\cong"]\\
				\Hom_{\Sh(\BBB')}((g \circ w)^s(h^{\#}_A),F) \arrow[r,"- \: \circ (\alpha^s)_{h^{\#}_A}"] &\Hom_{\Sh(\BBB')}((f\circ w)^s(h^{\#}_A), F)
			\end{tikzcd}
		\end{equation}
		is a commutative diagram. Now observe that:
		\begin{equation*}
		\begin{aligned}
		v \circ f  &= \#_{\AAA} \circ Y_{\AAA} \circ f = f^s \circ \#_{\BBB} \circ Y_{\BBB},\\
		v \circ g &= \#_{\AAA} \circ Y_{\AAA} \circ g = g^s \circ \#_{\BBB} \circ Y_{\BBB}.
		\end{aligned}
		\end{equation*}
		But notice that, as $w$ is an LC morphism, $w^s: \Sh(\BBB') \lra \Sh(\BBB)$ is an equivalence with quasi-inverse given by $w_s:\Sh(\BBB) \lra \Sh(\BBB')$. Hence we have that
		\begin{equation}\label{definingbeta}
		\begin{aligned}
		v \circ f &\cong (f \circ w)^s (w_s \circ \#_{\BBB} \circ Y_{\BBB})\\
		v \circ g &\cong (g \circ w)^s (w_s \circ \#_{\BBB} \circ Y_{\BBB}).
		\end{aligned}
		\end{equation}
		We define $\beta$ as the composition
		\begin{equation}
		\begin{tikzcd}
		v \circ f \arrow[r,Rightarrow,"\cong"] &(f \circ w)^s \circ w_s \circ \#_{\BBB} \circ Y_{\BBB} \arrow[r,Rightarrow] &(g \circ w)^s \circ w_s \circ \#_{\BBB} \circ Y_{\BBB} \arrow[r,Rightarrow,"\cong"] &v \circ g
		\end{tikzcd}
		\end{equation}
		where the second 2-morphism is given by $\alpha^s \circ (w_s \circ \#_{\BBB} \circ Y_{\BBB})$.
		
		Let's now check that $v \circ \alpha = \beta \circ w$. First observe that have the following commutative diagram
		\begin{equation*}
		\begin{tikzcd}[column sep=50pt]
			h^{\#}_{(g \circ w)(B')} \arrow[r,"\beta_{w(B')}"] \arrow[d,"\cong"] &h^{\#}_{(f \circ w)(B')} \arrow[d,"\cong"]\\
			(g\circ w)^s \circ w_s \circ h^{\#}_{w(B')} \arrow[r,"(\alpha^s)_{w_s (h^{\#}_{w(B')})}"] \arrow[d,"\cong"] &(f\circ w)^s \circ w_s \circ h^{\#}_{w(B')} \arrow[d,"\cong"]\\
			(g \circ w)^s(h^{\#}_{B'}) \arrow[r,"(\alpha^s)_{h^{\#}_{B'}}"] \arrow[d,equal] &(f\circ w)^s(h^{\#}_{B'}) \arrow[d,equal]\\
			h^{\#}_{(g \circ w)(B')} \arrow[r,"h^{\#}_{\alpha_{B'}}"] &h^{\#}_{(f \circ w)(B')},
		\end{tikzcd}
		\end{equation*}
		where the last commutative square comes from the commutative diagram (\ref{alphasrepres}) above. It is easily seen  that the vertical compositions are the identity. Hence we have that $\beta_{w(B')} = h^{\#}_{\alpha_{B'}} = v(\alpha_{B'})$ for all $B' \in \BBB'$, which concludes the argument.
			
		We prove now (2). Assume $\alpha$ is an invertible $2$-morphism. Then so are $\alpha_s$ and $\alpha^s$, and hence $\alpha^s \circ (t_s \circ w)$ is also an invertible $2$-morphism. As $\beta$ is obtained from $\alpha^s \circ (w_s \circ v)$ via pre- and postcomposing (vertically) with invertible $2$-morphisms, we conclude the argument.
		
		Finally, we prove (3). Assume there exists another $v': \AAA \lra \AAA'$ in $\LC$ and another 2-morphism $\beta': v' \circ f \Rightarrow v'\circ g$ with $v' \circ \alpha = \beta' \circ w$. We have to prove that there exist morphisms $u,u': \AAA' \lra \AAA''$ such that $u \circ v \in \LC$ and $u' \circ v' \in \LC$, and an invertible $2$-morphism $\epsilon: u \circ v \Rightarrow u'\circ v'$ such that the following diagram
		\begin{equation}\label{LF4part3}
		\begin{tikzcd}
		u \circ v \circ f \arrow[r,"u \circ \beta"] \arrow[d,"\epsilon \circ f"] &u \circ v \circ g \arrow[d,"\epsilon \circ g"]\\
		u' \circ v' \circ f \arrow[r,"u' \circ \beta'"] &u' \circ v' \circ g 
		\end{tikzcd}
		\end{equation}
		commutes. 
		
		Consider the equivalence $v'^s \circ v_s: \Sh(\AAA') \lra \Sh(\AAA')$,
		whose quasi-inverse is given by $v^s \circ v'_s$. We consider the associated roof construction for $v^s \circ v'_s$:
		\begin{equation*}
		\begin{tikzcd}
		&\AAA'' \\
		\AAA' \arrow[ur,"u'"] &&\AAA', \arrow[ul,"u"'] 
		\end{tikzcd}
		\end{equation*}	
		with $\AAA'' = \{ h^{\#}_{A'}\}_{A' \in \AAA'} \cup \{ v^s \circ v'_s (h^{\#}_{A'}) \}_{A' \in \AAA'}$, $u = \#_{\AAA'} \circ Y_{\AAA'}$ and $u' = (v^s \circ v'_s) \circ \#_{\AAA'} \circ Y_{\AAA'}$.
		In particular, $u$ is in $\LC$ and as so was $v$ by assumption, so it follows from \Cref{CondLF2} that $u \circ v$ belongs to $\LC$.
		
		Now let's construct $\epsilon$. For each $A \in \AAA$ we have:
		\begin{equation*}
		u \circ v = \#_{\AAA'} \circ Y_{\AAA'} \circ v = v^s \circ \#_{\AAA} \circ Y_{\AAA} \cong (v^s \circ v'_s \circ v'^s) \circ \#_{\AAA} \circ Y_{\AAA} \cong (v^s \circ v'_s)(\#_{\AAA'} \circ Y_{\AAA'} \circ v') = u' \circ v'.
		\end{equation*}
		Let's denote this composition of invertible $2$-morphisms by $\epsilon: u \circ v \Rightarrow u' \circ v'$. Observe first that from \Cref{CondLF5} above, it follows that $u' \circ v'$ also belongs to $\LC$. To finish the argument, it remains to check that the diagram (\ref{LF4part3}) above is commutative. Evaluating in any object $B \in \BBB$ we obtain the following commutative diagram.
			\begin{equation} \label{bigcomdiag} 
				\begin{tikzcd}[row sep=small, column sep=60pt]
					u \circ v \circ f (B) \arrow[r,"u (\beta_B)"] \arrow[d,equals] &u \circ v \circ g(B) \arrow[d,equals]\\
					h^{\#}_{v \circ f(B)} \arrow[r,"h^{\#}_{\beta_B}"] \arrow[d,"\cong"] &h^{\#}_{v \circ g(B)} \arrow[d,"\cong"]\\
					(v \circ f \circ w)^s (w_s(h^{\#}_B)) \arrow[r,"(\beta \circ w)^s_{w_s (h^{\#}_B)}"] \arrow[d,equals] &(v \circ g \circ w)^s  (w_s(h^{\#}_B)) \arrow[d,equals]\\
					(v \circ f \circ w)^s (w_s (h^{\#}_B)) \arrow[r,"(v \circ \alpha)^s_{w_s h^{\#}_B}"] \arrow[d,"\cong"] &(v \circ f \circ w)^s (w_s (h^{\#}_B)) \arrow[d,"\cong"]\\
					(v^s \circ v'_s) \circ (v'^s \circ f \circ w)^s (w_s (h^{\#}_B)) \arrow[r,"(v^s \circ v'_s) \circ (v' \circ \alpha)^s_{w_s h^{\#}_B}"] \arrow[d,equals] &(v \circ v'_s) \circ (v'^s \circ f \circ w)^s (w_s (h^{\#}_B)) \arrow[d,equals]\\
					(v^s \circ v'_s) \circ (v'^s \circ f \circ w)^s (w_s (h^{\#}_B)) \arrow[r,"(v^s \circ v'_s) \circ (\beta' \circ w)^s_{w_s h^{\#}_B}"] \arrow[d,"\cong"] &(v^s \circ v'_s) \circ (v'^s \circ g \circ w)^s (w_s (h^{\#}_B)) \arrow[d,"\cong"]\\
					(v^s \circ v'_s) \circ (v'^s \circ f)^s (h^{\#}_B)\arrow[r,"(v^s \circ v'_s) \circ \beta'^s_{h^{\#}_B}"] \arrow[d,equal] &(v^s \circ v'_s) \circ (v'^s \circ g)^s (h^{\#}_B) \arrow[d,equal]\\
					(v^s \circ v'_s) (h^{\#}_{(v' \circ f)(B)})\arrow[r,"(v^s \circ v'_s) \circ h^{\#}_{\beta'_B}"] \arrow[d,equal] &(v^s \circ v'_s) (h^{\#}_{(v' \circ g)(B)}) \arrow[d,equal]\\
					u' \circ v' \circ f (B) \arrow[r,"u' (\beta'_B)"] &u' \circ v' \circ g(B)
				\end{tikzcd}
			\end{equation}
		Observe that the left vertical composition in the diagram equals $\epsilon_{f(B)}$ and the right vertical composition equals $\epsilon_{g(B)}$, which concludes our argument.	
	\end{proof}
\end{lemma}

Finally, we are in the position to prove the following.
\begin{proposition}
	$\LC$ admits a left calculus of fractions in $\site$.
	\begin{proof}
		The statement follows from \Cref{CondLF1,CondLF2,CondLF3,CondLF5,CondLF4} above.
	\end{proof}	
\end{proposition}
Hence we can localize $\site$ with respect to $\LC$ and obtain the bilocalization $\site[\LC^{-1}]$.

\section{The 2-category of Grothendieck categories as a bilocalization of the 2-category of sites}\label{sectionlocalizationallsites}
In this section we prove the main result of the paper.

\begin{theorem}\label{mainresultbis}
There exists a pseudofunctor
$$\Phi: \site \lra \groth$$
which sends LC morphisms to equivalences in $\groth$, such that the pseudofunctor
$$\tilde{\Phi}: \site[\LC^{-1}] \lra \groth$$
induced by $\Phi$ via the universal property of the bicategory of fractions is an equivalence of bicategories.
\end{theorem}

Let $\ccc$ be a bicategory and $\mathsf{W}$ a class of $1$-morphisms in $\ccc$ that admits a calculus of left fractions. Given a bicategory $\ddd$ and a pseudofunctor $\Phi: \ccc \lra \ddd$ sending $1$-morphisms that belong to $\mathsf{W}$ to equivalences in $\ddd$, we have that $\Phi$ induces a pseudofunctor $\tilde{\Phi}: \ccc[\mathsf{W}^{-1}] \lra \ddd$ by the universal property of bilocalizations.
A characterization of the pseudofunctors $\Phi$ such that $\tilde{\Phi}$ is an equivalence of bicategories (in the case of a right bicategory of fractions) is provided in \cite{tommasini3} . The characterization makes use of the \emph{right saturation} of a class of morphisms introduced in \cite{tommasini2}. We formulate below an analogue for a left calculus of fractions.

\begin{definition}
	Let $\mathsf{W}$ be a class of $1$-morphisms in the bicategory $\ccc$. The \emph{(left) saturation} $\mathsf{W}_{\mathrm{sat}}$ of $\mathsf{W}$ is the class of all $1$-morphisms $f : A \lra B$ in $\ccc$, such that there exists a pair of objects $C,D \in \ccc$ and a pair of morphisms $g : B \lra C$ and $h : C \lra D$, such that both $g \circ f$ and $h \circ g$ belong to $\mathsf{W}$.
	
	We say that $\mathsf{W}$ is \emph{(left) saturated} if $\mathsf{W} = \mathsf{W}_{\mathrm{sat}}$.
\end{definition}

In analogy to \cite[Rem 2.3]{tommasini2} in the case of right saturation, we have the following statement for the left saturation.
\begin{proposition}\label{leftsaturation}
	If the class of morphisms $\mathsf{W}$ admits a left calculus of fractions, then $\mathsf{W} \subseteq \mathsf{W}_{\mathrm{sat}}$.
\end{proposition}

\begin{proposition}[{\cite[Thm 0.4]{tommasini3}}]\label{mainresultbilocalization}
	Let $\ccc$ be a bicategory and $\mathsf{W}$ a class of $1$-morphisms in $\ccc$ that admits a calculus of left fractions. Given a bicategory $\ddd$ and a pseudofunctor $\Phi: \ccc \lra \ddd$ sending $1$-morphisms that belong to $\mathsf{W}$ to equivalences in $\ddd$, we have that $\Phi$ induces an equivalence of bicategories $\tilde{\Phi}: \ccc[\mathsf{W}^{-1}] \overset{\cong}{\lra} \ddd$ if and only if:
	\begin{itemize}
		\item[\emph{\textbf{B1}}] $\Phi$ is essentially surjective on objects;
		\item[\emph{\textbf{B2}}]  Given objects $C_1,C_2 \in \ccc$ and an equivalence $e: \Phi(C_2) \overset{\cong}{\lra} \Phi(C_1)$, there exits an object $C_3 \in \ccc$, a pair of morphisms $w_1: C_1 \lra C_3$ in $\mathsf{W}$ and $w_2: C_2 \lra C_3$ in $\mathsf{W}_{\mathrm{sat}}$, an equivalence $e':\Phi(C_3) \lra \Phi(C_1)$ in $\ddd$ and a pair of invertible $2$-morphisms $\delta_1,\delta_2$ as follows
		\begin{equation}
		\begin{tikzcd}
		\Phi(C_2) \arrow[dr,"\Phi(w_2)"'] \arrow[drr, bend left, "e"] \arrow[drr, phantom, description, "\Downarrow \delta_2", pos=0.6,yshift=1ex]\\
		&\Phi(C_3) \arrow[r, "e'"]  &\Phi(C_1).\\
		\Phi(C_1)\arrow[ur,"\Phi(w_1)"] \arrow[urr, bend right,"\id_{\Phi(C_1)}"'] \arrow[urr, phantom, description,  "\Downarrow \delta_1", pos=0.6, yshift=-1ex]
		\end{tikzcd}
		\end{equation} 
		\item[\emph{\textbf{B3}}]  Given objects $C \in \ccc$ and $D \in \ddd$ and a morphism $f: \Phi(C) \lra D$, there exists an object $C' \in \ccc$, a morphism $g: C \lra C'$ in $\ccc$, an equivalence $e:\Phi(C') \overset{\cong}{\lra} D$ in $\ddd$ and an invertible $2$-morphism $\alpha: f \Rightarrow e \circ \Phi(g)$.
		\item[\emph{\textbf{B4}}]  Given objects $C, C' \in \ccc$, two $1$-morphisms $f_1, f_2: C \lra C'$ in $\ccc$ and two $2$-morphisms $\gamma_1,\gamma_2: f_1 \Rightarrow f_2$ such that $\Phi(\gamma_1) = \Phi(\gamma_2)$, there exits an object $C'' \in \ccc$ and a $1$-morphism $w: C' \lra C''$ in $\mathsf{W}$ such that $w \circ \gamma_1 = w \circ \gamma_2$.
		\item[\emph{\textbf{B5}}]  Given objects $C,C' \in \ccc$, a pair of morphisms $f_1,f_2: C \lra C'$ and a $2$-morphism $\alpha: \Phi(f_1) \Rightarrow \Phi(f_2)$, then, there is an object $C'' \in \ccc$, a morphism $w: C' \lra C''$ in $\mathsf{W}$ and a $2$-morphism $\beta: w \circ f_1 \Rightarrow w \circ f_2$ such that $\Phi(w) \circ \alpha = \psi^{\Phi}_{w,f_2} \bullet \Phi(\beta) \bullet (\psi^{\Phi}_{w,f_1})^{-1} $, where $\psi^{\Phi}$ denotes the associator of the pseudofunctor $\Phi$.
	\end{itemize}
\end{proposition} 

\begin{remark}
	We need to use this set of necessary and sufficient conditions from \cite{tommasini3} as the set of sufficient conditions provided by \cite[Prop 24]{pronk} is not satisfied in our case.
\end{remark}

We are now in the position to prove \Cref{mainresultbis}.
\begin{proof}
Recall from \S\ref{section2catssitesgroth} that we have a pseudofunctor $\Phi: \site \lra \groth$ that sends LC morphisms to equivalences. Hence by the universal property of the bilocalization, we have an induced pseudofunctor
\begin{equation}
\tilde{\Phi}: \site[\LC^{-1}] \lra \groth.
\end{equation}
Consequently, if the pseudofunctor $\Phi: \site \lra \groth$ satisfies properties \emph{\textbf{B1}} to \emph{\textbf{B5}} from \Cref{mainresultbilocalization}, we conclude the argument. This is done in \Cref{propertiesB} below.
\end{proof}

\begin{lemma}\label{propertiesB}
	The pseudofunctor $\Phi: \site \lra \groth$ satisfies properties \emph{\textbf{B1}} to \emph{\textbf{B5}} in \Cref{mainresultbilocalization} above.
	\begin{proof}
		We know every Grothendieck category can be realised as a category of sheaves on a site (see \S \ref{seclinearsites}), hence $\Phi$ is essentially surjective on objects, which proves \textbf{B1}.
		
		We now prove \textbf{B2}. Consider two sites $\AAA_1,\AAA_2$ and an equivalence $e: \Sh(\AAA_2) \overset{\cong}{\lra} \Sh(\AAA_1)$ in $\groth$. We apply the roof theorem to the functor $e$. Let $\AAA_3$ be the site with objects $$\Obj(\AAA_3) =\{h^{\#}_{A_1}\}_{A_1 \in \AAA_1} \cup \{e(h^{\#}_{A_2})\}_{A_2 \in \AAA_2} \subseteq \Sh(\AAA_1),$$ and the topology induced by the canonical topology in $\Sh(\AAA_1)$. We have the roof construction
		\begin{equation*}
		\begin{tikzcd}
		&\AAA_3 \\
		\AAA_2 \arrow[ur,"w_2"] &&\AAA_1, \arrow[ul,"w_1"'] 
		\end{tikzcd}
		\end{equation*}
		where $w_1 = \#_{\AAA_1} \circ Y_{\AAA_1}$ and $w_2= e \circ \#_{\AAA_2} \circ Y_{\AAA_2}$. By the roof theorem, $w_1$ is an LC morphism. On the other hand, 
		we have that
		$$e \cong \widetilde{(w_1)}^* \circ (w_2)^s \cong (w_1)_s \circ (w_2)^s,$$
		where the first step is given by the roof theorem and the second by \Cref{cocontandcontLC}. Observe that $(w_1)_s$ is an equivalence because $w_1$ is an LC morphism. Then, as $e$ is an equivalence by hypothesis, $(w_2)^s$ is also and equivalence and hence so is $(w_2)_s$. In summary, we have a morphism $w_2: \AAA_2 \lra \AAA_3$ with $\AAA_3$ subcanonical and such that $(w_2)_s$ is an equivalence. Then, it follows from \cite[Cor 4.5]{lowenGP} that $w_2$ is an LC morphism. 
	 	Hence, in particular, as $\LC \subseteq \LC_{\mathrm{sat}}$ by \Cref{leftsaturation}, we have that $w_2 \in \LC_{\mathrm{sat}}$.
		Consider now the equivalence $$e' : \Sh(\AAA_3) \overset{\cong}{\lra} \Sh(\AAA_1),$$ given by $e'= \widetilde{(w_1)}^*$. Then we can then choose $\delta_2$ to be the isomorphism:
		$$e \cong e' \circ (w_2)^s = e' \circ \Phi(w_2)$$ 
		given by the roof decomposition of $e$. On the other hand, we have that
		$$e' \circ \Phi(w_1) = \widetilde{(w_1)}^* \circ (w_1)^s \cong (w_1)_s \circ (w_1)^s \cong \id_{\Sh(\AAA_1)},$$ 
		where the second step follows from \Cref{cocontandcontLC}, and the last from the fact that $w_1$ is LC, and hence $(w_1)^s$ is an equivalence with quasi-inverse given by $(w_1)_s$. Then, we can denote by $\delta_1$ the horizontal composition of this chain of invertible $2$-morphisms, which concludes the argument.
		
		We now proceed to prove \textbf{B3}. Fix a site $\BBB$, a Grothendieck category $\aaa$ and a $1$-morphism $f:\Sh(\BBB) \lra \aaa$ in $\groth$. We choose a site $\AAA$ such that we have an equivalence $e': \aaa \overset{\cong}{\lra} \Sh(\AAA)$. Consider the morphism $e' \circ f:\Sh(\BBB) \lra \Sh(\AAA)$ in $\groth$ and its associated roof decomposition
		\begin{equation*}
		\begin{tikzcd}
		&\CCC\\
		\BBB \arrow[ur,"g"] &&\AAA. \arrow[ul,"w",swap]
		\end{tikzcd}
		\end{equation*}
	 	We then have an equivalence $\tilde{w}_* \circ e':\aaa \overset{\cong}{\lra} \Sh(\AAA) \overset{\cong}{\lra} \Sh(\CCC)$. Choose $e''$ a quasi-inverse of $e'$ and consider $e = e'' \circ  \tilde{w}^*: \Sh(\CCC) \overset{\cong}{\lra} \aaa$ which is a quasi-inverse of $\tilde{w}^* \circ e'$. Then, by the roof theorem, we have that $e' \circ f \cong \tilde{w}^* \circ g^s$ and hence, by postcomposing with $e''$ on both the left and the right hand side, we have an invertible $2$-morphism $f \cong e \circ g^s = e \circ \Phi(g)$, which finishes the argument.
		
		We now prove \textbf{B4}. Fix two sites $\AAA, \BBB$, two continuous morphisms $f_1,f_2: \AAA \lra \BBB$ and a pair of 2-morphisms $\gamma_1,\gamma_2:f_1 \Rightarrow f_2$ such that $$(\gamma_1)^s = (\gamma_2)^s: (f_1)^s \Rightarrow (f_2)^s.$$ Consider $\CCC$ the site with objects $\{h^{\#}_{B}\}_{B \in \BBB}$ with the topology induced by the canonical topology in $\Sh(\BBB)$ and $w = \#_{\BBB} \circ Y_{\BBB}: \BBB \lra \CCC$ the corresponding LC morphism. Observe that $((\gamma_1)^s)_{h^{\#}_A} = ((\gamma_2)^s)_{h^{\#}_A}$ for all $A \in \AAA$. This implies, applying the commutative diagram \eqref{alphasrepres}, that $h^{\#}_{(\gamma_1)_A} = h^{\#}_{(\gamma_2)_A}$ for all $A \in \AAA$. Hence, we have that $(w \circ \gamma_1)_A = (w \circ \gamma_2)_A$ for all $A \in \AAA$ and natural in $A$, which concludes the argument. 
		
		Finally, we prove property \textbf{B5}. Consider two sites $\AAA, \BBB$, two continuous morphisms $f_1,f_2: \AAA \lra \BBB$ and a $2$-morphism $\alpha: f_1^s \Rightarrow f_2^s$. Let $\CCC$ the site with objects $\{h^{\#}_{B}\}_{B \in \BBB} \subseteq \Sh(\BBB)$ and the topology induced by the canonical topology in $\Sh(\BBB)$ and the LC morphism $w = \#_{\BBB} \circ Y_{\BBB}: \BBB \lra \CCC$. Take $\beta: w\circ f_1 \Rightarrow w \circ f_2$ the $2$-morphism  given by $$w \circ f_1 = \#_{\BBB} \circ Y_{\BBB} \circ f_1 = f_1^s \circ \#_{\AAA} \circ Y_{\AAA} \overset{\alpha \circ \#_{\AAA} \circ Y_{\AAA}}{\Longrightarrow} f_2^s \circ \#_{\AAA} \circ Y_{\AAA} = \#_{\BBB} \circ Y_{\BBB} \circ f_2 = w \circ f_2.$$
		Then, we have that 
		\begin{equation*}
		\begin{tikzcd}[column sep= 80pt]
		w^s \circ f_1^s \arrow[r,Rightarrow,"w^s \circ \alpha"]  &w^s \circ f_2^s \arrow[d, Rightarrow, "(\psi^{S}_{w,f_2})^{-1}", "\cong"'] \\
		(w \circ f_1)^s \arrow[r,Rightarrow,"(\psi^{S}_{w,f_2})^{-1} \bullet (w^s \circ \alpha)\bullet \psi^{S}_{w,f_1}"] \arrow[u, Rightarrow, "\psi^{S}_{w,f_1}", "\cong"'] \arrow[d,equal] &	(w \circ f_2)^s \arrow[d,equal]  \\
		(f_1^s \circ \#_{\AAA} \circ Y_{\AAA})^s \arrow[r,Rightarrow,"(\alpha \circ\#_{\AAA} \circ Y_{\AAA})^s "]&(f_2^s \circ \#_{\AAA} \circ Y_{\AAA})^s
		\end{tikzcd}
		\end{equation*}
		is a commutative diagram of $2$-morphisms.
		But observe that the composition $$(w \circ f_1)^s = (f_1^s \circ \#_{\AAA} \circ Y_{\AAA})^s \overset{(\alpha \circ\#_{\AAA} \circ Y_{\AAA})^s}{\Longrightarrow} (f_2 \circ \#_{\AAA} \circ Y_{\AAA})^s = (w \circ f_2)^s$$ is just $\beta^s$ by definition, hence
		$$w^s \circ \alpha = \psi^{S}_{w,f_2} \bullet \beta^s \bullet (\psi^{S}_{w,f_1})^{-1},$$
		which concludes the argument.
	\end{proof}
\end{lemma}

\section{Monoidal bilocalization}\label{tensorproductvialocalization}
Let $\ccc$ be a category and $\mathsf{W}$ a class of morphisms which admits a calculus of fractions in the sense of Gabriel-Zisman \cite{gabriel-zisman}. Then, as it is proven in \cite{day}, if $\ccc$ has a symmetrical monoidal structure such that $\mathsf{W}$ is closed under tensoring, then the localization $\ccc[\mathsf{W}^{-1}]$ has a monoidal structure such that the localization functor $\ccc \lra \ccc[\mathsf{W}^{-1}]$ is a monoidal functor. This is what in \cite{day} is referred to as \emph{monoidal localization}.

It is reasonable to believe that an analogous result for monoidal bicategories and bilocalizations holds true, and we plan to return to this topic in the future. In this section, we briefly sketch a possible application of the main result of this paper given that we have a satisfactory theory of \emph{monoidal bilocalization} available.

In \cite{lowen-ramos-shoikhet}, we define a tensor product of linear sites, which is seen to define a symmetric monoidal structure on the bicategory $\site$
$$\boxtimes: \site \times \site \lra \site: ((\AAA, \ttt_{\AAA}), (\BBB, \ttt_{\BBB})) \longmapsto (\AAA \otimes \BBB, \ttt_{\AAA} \boxtimes \ttt_{\BBB}).$$
We further showed in loc. cit. that the class of LC morphisms is closed under $\boxtimes$, hence we obtain an induced bi-pseudofunctor
$$\tilde{\boxtimes}: \site[\LC^{-1}] \times \site[\LC^{-1}] = (\site \times \site)[(\LC \times \LC)^{-1}] \lra \site[\LC^{-1}]$$
which a general theory of monoidal bilocalization would yield to define a monoidal structure in the bicategorical sense. This structure could then be transferred to the equivalent bicategory $\groth$ (in a non-canonical way).

Note that in \cite{lowen-ramos-shoikhet}, we use the roof construction in order to obtain a well-defined (up to equivalence) tensor product of Grothendieck categories, at least on the level of the categories, and we show that this tensor product is a special case of the tensor product of locally presentable categories which is known to be bi-functorial and monoidal. Further, we provide an alternative approach to these issues in \cite{ramosthesis}, making use of canonical bicolimit presentations of Grothendieck categories.

\def\cprime{$'$} \def\cprime{$'$}
\providecommand{\bysame}{\leavevmode\hbox to3em{\hrulefill}\thinspace}
\providecommand{\MR}{\relax\ifhmode\unskip\space\fi MR }
\providecommand{\MRhref}[2]{%
	\href{http://www.ams.org/mathscinet-getitem?mr=#1}{#2}
}
\providecommand{\href}[2]{#2}

\end{document}